\documentclass[11pt, a4]{article}

\usepackage{amsmath}
\usepackage{amssymb}
\usepackage{amscd}

\newtheorem{theorem}{Theorem}[section]

\newtheorem{definition}[theorem]{Definition}

\title{\bf A sound interpretation of Le\'{s}niewski's epsilon in modal logic \bf KTB}
\author{Takao Inou\'{e}\\ \\This is the second version, October 2, 2021.}
\date{}
\begin{document}

\maketitle

\begin{abstract}
  In this paper, we shall show that the following translation $I^M$ from the propositional 
fragment $\bf L_1$ of Le\'{s}niewski's ontology to modal logic 
$\bf KTB$ is sound: for any formula $\phi$ and $\psi$ of $\bf L_1$, it is defined as 

\smallskip 

(M1) $I^M(\phi \vee \psi)$ = $I^M(\phi) \vee I^M(\psi),$

(M2) $I^M(\neg \phi)$ = $\neg I^M(\phi),$

(M3) $I^M(\epsilon ab)$ = $\Diamond p_a \supset p_a . \wedge . \Box p_a \supset \Box p_b . 
\wedge . \Diamond p_b \supset p_a,$

\smallskip 

\noindent where $p_a$ and $p_b$ are propositional variables corresponding to the name variables $a$ and $b$, 
respectively. In the last section, we shall give some open problems and my conjectures.
\end{abstract}

\smallskip

\normalsize 

\section{Introduction and $I^M$}

Inou\'{e} \cite{inoue1} initiated a study of interpretations of Le\'{s}niewski's epsion $\epsilon$ in the modal logic $\bf K$ and its certain extensions. That is, Ishimoto's propositional fragment $\bf L_1$ (Ishimoto \cite{ishi}) of Le\'{s}niewski's ontology $\bf L$ (refer to Urbaniak \cite{urbaniak-book}) is partially embedded 
in $\bf K$ and in its certain extensions, respectively, by the following translation $I$ from $\bf L_1$ to them: 
for any formula $\phi$ and $\psi$ of $\bf L_1$, it is defined as 

\smallskip 

(I1) $I(\phi \vee \psi)$ = $I(\phi) \vee I(\psi)$, 

(I2) $I(\neg \phi)$ = $\neg I(\phi)$, 

(I3) $I(\epsilon ab)$ = $p_a \wedge \Box (p_a \equiv p_b)$,

\smallskip 

\noindent where $p_a$ and $p_b$ are propositional variables corresponding to the name variables $a$ and $b$, respectively. 
Here, ``$\bf L_1$ is partially embedded in $\bf K$ by $I$" means that for any formula $\phi$ 
of a certain decidable nonempty set of formulas of $\bf L_1$ (i.e. decent formulas (see $\S 3$ of Inou\'{e} \cite{inoue3})), $\phi$ is 
a theorem of $\bf L_1$ if and only if $I(\phi)$ is a theorem of $\bf K$. Note that $I$ is sound. 
The paper \cite{inoue3} also proposed similar partial interpretations of Le\'{s}niewski's epsilon 
in certain von Wright-type deontic logics, that is, ten Smiley-Hanson systems of monadic 
deontic logic and in provability logic $\bf GL$, respectively. (See \AA qvist \cite{aq} and 
Boolos \cite{boolos} for those logics.)

The interpretation $I$ is however not faithful. A counterexample for the faithfulness is, for example, 
$\epsilon ac \wedge \epsilon bc. \supset . \epsilon ab \vee \epsilon cc$ (for the details, see \cite{inoue3}). Blass \cite{blass} 
gave a modification of the interpretation and showed that his interpretation $T$ is faithful, using Kripke models. 
Inou\'{e} \cite{inoue-Blass} called the translation \it Blass translation \rm (for short, \it $B$-translation\rm ) or 
\it Blass interpretation \rm (for short, \it B-interpretation\rm ). The translation $B$ from $\bf L_1$ to $\bf K$ is 
defined as follows: for any formula 
$\phi$ and $\psi$ of $\bf L_1$, 

\smallskip 

(B1) $B(\phi \vee \psi)$ = $B(\phi) \vee B(\psi)$, 

(B2) $B(\neg \phi)$ = $\neg B(\phi)$, 

(B3) $B(\epsilon ab)$ = $p_a \wedge \Box (p_a \supset p_b) \wedge . p_b \supset \Box (p_b \supset p_a)$,

\smallskip 

\noindent where $p_a$ and $p_b$ are propositional variables corresponding to the name variables $a$ and $b$, respectively. 
Inou\'{e} \cite{inoue-Blass} extended Blass's faithfulness result for many normal modal logics, provability logic and von Wright-type 
dentic logics including $\bf K4$, $\bf KD$, $\bf KB$, $\bf KD4$, etc, $\bf GL$ and ten Smiley-Hanson systems of monadic deontic logic, 
using model constructions based on Hintikka formula.

In this paper, we first propose a translation $I^M$ from $\bf L_1$ in modal logic $\bf KTB$, which will be specified in \S2.

\begin{definition}
A translation $I^M$ of Le\'{s}niewski's propositional ontology 
$\bf L_1$ in modal logic $\bf KTB$ is defined as follows: for any formula $\phi$ and $\psi$ of $\bf L_1$, 

\smallskip \rm 

(M1) $I^M(\phi \vee \psi)$ = $I^M(\phi) \vee I^M(\psi),$

(M2) $I^M(\neg \phi)$ = $\neg I^M(\phi),$

(M3) $I^M(\epsilon ab)$ = 
$\Diamond p_a \supset p_a. \wedge . \Box p_a \supset \Box p_b . \wedge . \Diamond p_b \supset p_a,$

\smallskip \it 

\noindent where $p_a$ and $p_b$ are propositional variables corresponding to the name variables $a$ and $b$, respectively. 

\end{definition}

We call $I^M$ to be \it M-translation \rm or \it M-interpretation\rm . 

In the following \S2, we shall collect the basic preliminaries for this paper.  In \S3, using proof theory, we 
shall show that $I^M$ is sound, as the main theorem of this paper. In S4, we shall give some comments including 
some open problems and my conjectures.

%%%%%%%%%%%%%%%%%%%%%%%%%%%%%%%

\section{Propositional ontology $\bf L_1$ and modal logic KTB}

Let us recall a formulation of $\bf L_1$, which was introduced in \cite{ishi}. The Hilbert-style system of it, denoted again by 
$\bf L_1$, consists of the following axiom-schemata with a formulation of classical propositional logic $\bf CP$ as its axiomatic basis:

\smallskip 

(Ax1) $\enspace$ $\epsilon ab$ $\supset$ $\epsilon aa$,

\smallskip 

(Ax2) $\enspace$ $\epsilon ab$ $\wedge$ $\epsilon bc$. $\supset$ $\epsilon ac$,

\smallskip 

(Ax3) $\enspace$ $\epsilon ab$ $\wedge$ $\epsilon bc$. $\supset$ $\epsilon ba$,

\smallskip 

\noindent where we note that every atomic formula of $\bf L_1$ is of the form $\epsilon ab$ for some name variables $a$ and $b$ 
and a possible intuitive interpretation of $\epsilon ab$ is `the $a$ is $b$'. We note that (Ax1), (Ax2) and (Ax3) are theorems of 
Le\'{s}niewski's ontology (see S\l upecki \cite{slu}). 

The modal logic $\bf K$ is the smallest logic which contains all instances of classical tautology and all formulas of the forms 
$\Box (\phi \supset \psi) \supset . \Box \phi \supset \Box \psi$ being closed under modus ponens and the rule of necessitation 
(for $\bf K$ and basics for modal logic, see Bull and Segerberg \cite{bulseg}, Chagrov and Zakharyaschev \cite{cz}, Fitting 
\cite{fitting}, Hughes and Cresswell \cite{hc} and so on). 

We recall the naming of modal logics as follows (refer to e.g. Poggiolesi \cite{pog} and Ono \cite{ono1}, 
also see Bull and Segerberg \cite{bulseg}):

\smallskip

%$\bf KD$: $\bf K$ + $\Box \phi \supset \Diamond \phi$ ($\bf D$, serial relation)

$\bf KT$: $\bf K$ + $\Box \phi \supset \phi$ ($\bf T$, reflexive relation)

$\bf KB$: $\bf K$ + $\phi \supset \Box \Diamond \phi$ ($\bf B$, symmetric relation)

$\bf KTB$: $\bf KT$ + $\bf B$ (reflexive and symmetric relation).

%$\bf KTBD$: $\bf KT$ + $\bf B$ + $\bf D$

%\smallskip

%\noindent Note that $\bf KT$ is properly stronger than 
%$\bf KD$, since a reflexive  relation  is a serial one. So, $\bf KTBD$ is the same logic as $\bf KTB$.

%%%%%%%%%%%%%%%%%%%%%%%%%%%%%%%%%%%%%%

\section{The soundness of $I^M$}

\begin{theorem} \rm (Soundness) 
\it For any formula $\phi$ of $\bf L_1$, we have  
$$\vdash_{\bf L_1} \phi \enspace \Rightarrow \enspace \vdash_{\bf KTB} I^M(\phi).$$
\end{theorem}

\noindent \sc Proof. \rm  Let $\phi$ be a formula of $\bf L_1$.  We shall prove the meta-implication by induction on derivation.

\noindent \sc Basis. \rm 

\smallskip

\noindent (Case 1) We shall first treat the case for (Ax1). Let $a$ and $b$ be name variables. 
Then we have the following inferences in $\bf KTB$:

\smallskip

$(*)$ $I^M(\epsilon ab)$ (Assumption)

(1.1) $\Diamond p_a \supset p_a$ from $(*)$ and Definition 1.1) \dag 

(1.2) $\Box p_a \supset \Box p_a$ (true in $\bf K$) \dag 

(1.3) $\Diamond p_a \supset p_a. \wedge . \Box p_a \supset \Box p_a . \wedge . \Diamond p_a \supset p_a$ (from (1.1) and (1.2))

(1.4) $I^M(\epsilon aa)$ (from (1.3) and Definition 1.1)

(1.5) $I^M(\epsilon ab \supset \epsilon aa)$ (from $(*)$, (1.4) and Definition 1.1).

\smallskip

\noindent (Case 2) Next we shall deal with the case of (Ax2). Let $a$, $b$ and $c$ be name variables. 
Then we have the following inferences in $\bf KTB$:

\smallskip

$(**)$ $I^M(\epsilon ab \wedge \epsilon bc)$ (Assumption)

(2.1) $I^M(\epsilon ab)$ (from $(**)$ and Definition 1.1)

(2.2) $I^M(\epsilon bc)$ (from $(**)$ and Definition 1.1)

(2.3) $\Diamond p_a \supset p_a. \wedge . \Box p_a \supset \Box p_b . \wedge . \Diamond p_b \supset p_a$ (from (2.1) and Def 1.1)

(2.4) $\Diamond p_b \supset p_b. \wedge . \Box p_b \supset \Box p_c . \wedge . \Diamond p_c \supset p_b$ (from (2.2) and Def 1.1)

(2.5) $\Diamond p_a \supset p_a$ (from (2.3)) \dag 

(2.6) $\Box p_a \supset \Box p_b$ (from (2.3)) 

(2.7) $\Box p_b \supset \Box p_c$ (from (2.4))

(2.8) $\Box p_a \supset \Box p_c$ (from (2.6) and (2.7)) \dag 

(2.9) $\Diamond p_b \supset p_a$ (from (2.3))

(2.10) $\Box(\Diamond p_b \supset p_a)$ (from (2.9) and the rule of necessitation)

(2.11) $\Box \Diamond p_b \supset \Box p_a$ (from (2.10) with a true inference in $\bf K$)

(2.12) $\Box p_a \supset p_a$ (true in $\bf KT$)

(2.13) $\Box \Diamond p_b \supset p_a$ (from (2.11) and (2.12))

(2.14) $p_b \supset \Box \Diamond p_b$ (true in $\bf KB$)

(2.15) $\Diamond p_c \supset p_b$ (from (2.4))

(2.16) $\Diamond p_c \supset p_a$ (from (2.13) and (2.14) and (2.15)) \dag

(2.17) $\Diamond p_a \supset p_a. \wedge . \Box p_a \supset \Box p_c . \wedge . \Diamond p_c \supset p_a$ (from (2.5), (2.8) and (2.16))

(2.18) $I^M(\epsilon ac)$ (from (2.17) and Definition 1.1)

(2.19) $I^M(\epsilon ab \wedge \epsilon bc. \supset \epsilon ac)$ (from $(**)$, (2.18) and Definition 1.1).

\smallskip

\noindent (Case 3) Lastly we shall proceed to the case of (Ax3). Let $a$, $b$ and $c$ be name variables. 
Then we also have the following inferences in $\bf KTB$:

$(***)$ $I^M(\epsilon ab \wedge \epsilon bc)$ (Assumption)

(3.1) $I^M(\epsilon ab)$ (from $(***)$ and Definition 1.1)

(3.2) $I^M(\epsilon bc)$ (from $(***)$ and Definition 1.1)

(3.3) $\Diamond p_a \supset p_a. \wedge . \Box p_a \supset \Box p_b . \wedge . \Diamond p_b \supset p_a$ (from (3.1) and Def 1.1)

(3.4) $\Diamond p_b \supset p_b. \wedge . \Box p_b \supset \Box p_c . \wedge . \Diamond p_c \supset p_b$ (from (3.2) and Def 1.1)

(3.5) $\Diamond p_b \supset p_b$ (from (3.4)) \dag 

(3.6) $\Diamond p_b \supset p_a$ (from (3.3))

(3.7) $\Box (\Diamond p_b \supset p_a)$ (from (3.6) and the rule of necessitation)

(3.8) $\Box \Diamond p_b \supset \Box p_a$ (from (3.7) with a true inference in $\bf K$)

(3.9) $p_b \supset \Box \Diamond p_b$ (true in $\bf KB$)

(3.10) $\Box p_b \supset p_b$ (true in $\bf KT$)

(3.11) $\Box p_b \supset \Box p_a$ (from (3.8) and (3.9) and (3.10)) \dag 

(3.12) $\Diamond p_a \supset p_a$ (from (3.3))

(3.13) $p_a \supset \Box \Diamond p_a$ (true in $\bf KB$)

(3.14) $\Diamond p_a \supset \Box \Diamond p_a$ (from (3.12) and (3.13))

(3.15) $\Box (\Diamond p_a \supset p_a)$ (from (3.12) and the rule of necessitation)

(3.16) $\Box \Diamond p_a \supset \Box p_a$ (from (3.15) with a true inference in $\bf K$)

(3.17) $\Diamond p_a \supset \Box p_a$  (from (3.14) and (3.16))

(3.18) $\Box p_a \supset \Box p_b$ (from (3.3))

(3.19) $\Diamond p_a \supset \Box p_b$ (from (3.17) and (3.18))

(3.20) $\Box p_b \supset p_b$ (true in $\bf KT$)

(3.21) $\Diamond p_a \supset p_b$ (from (3.19) and (3.20)) \dag 

(3.22) $\Diamond p_b \supset p_b. \wedge . \Box p_b \supset \Box p_a . \wedge . \Diamond p_a \supset p_b$ 

$ \enspace \enspace \enspace \enspace \enspace \enspace$ (from (3.5), (3.11) and (3.21))

(3.23) $I^M(\epsilon ba)$ (from (3.22) and Definition 1.1)

(3.24) $I^M(\epsilon ab \wedge \epsilon bc. \supset \epsilon ba)$ (from $(***)$, (3.23) and Definition 1.1).

\smallskip

\noindent \sc Induction Steps. \rm 
The induction step is easily 
dealt with.  Suppose that $\phi$ and $\phi \supset \psi$ are theorems of $\bf L_1$. By induction hypthesis, $I^M(\phi)$ and 
$I^M(\phi \supset \psi)$ ($\leftrightarrow I^M(\phi) \supset I^M(\psi)$) are theorems of $\bf KTB$.  By modus ponens, we obtain 
$\vdash_{\bf KTB} I^M(\psi)$.  Thus this completes the proof the theorem. $\blacksquare$

%%%%%%%%%%%%%%%%%%%%%%%%%%%%%

\section{Open problems and conjectures}

In this last section, we shall present several open problems and my conjectures.

\smallskip

\bf Open problem 1: \rm Is  $I^M$ faithful?

\bf Open problem 2: \rm Find the set of other translations and modal logics in which $\bf L_1$ is 
embedded. I think that there seems to be many possibilities.

\bf Open problem 3: \rm Can $\bf L_1$ be embedded in $\bf S4.2$? (See e.g. Hamkins and L\"{o}we \cite{ham-lo}.)

\bf Open problem 4: \rm Can $\bf L_1$ be embedded in Grzegorczyk's modal Logic?
 (See e.g. Savateev and Shamkanov \cite{savateev2021})

\smallskip

\noindent My conjectures are the following.

\smallskip

\bf Conjecture 1: \rm  $I^M$ is faithful.

\bf Conjecture 2: \rm  It seems that $\bf L_1$ can not be embedded in intuitionistic propositional logic.

\bf Conjecture 3: \rm  It seems that $\bf L_1$ can well be embedded in intuitionistic modal propositional logic. 

\bf Conjecture 4: \rm $I^M$ is an embedding of $\bf L_1$ in Grzegorczyk's modal Logic.

\bigskip 

\noindent Takao Inou\'{e}

\noindent Department of Medical Molecular Informatics

\noindent Meiji Pharmaceutical University

\noindent Tokyo, Japan

\medskip 

\noindent Graduate School of Science and Engineering

\noindent Hosei University

\noindent Tokyo, Japan

\medskip 

\noindent Department of Applied Informatics

\noindent Faculty of Science and Engineering

\noindent Hosei University

\noindent Tokyo, Japan

\medskip

\noindent ta-inoue@my-pharm.ac.jp

\noindent takao.inoue.22@hosei.ac.jp

\noindent takaoapple@gmail.com


\begin{thebibliography}{99}

\bibitem{aq} L.~\AA qvist, {\em Deontic logic}, pp. 605--714, in \cite{gab} and 
also pp. 147--264 (as a revised version) in \cite{gab2}. 

\bibitem{blass} A.~Blass, {\em A faithful modal interpretation of propositional ontology}, 
\bf Mathematica Japonica, \rm 40 (1994), pp. 217--223.

\bibitem{boolos} G.~Boolos, \bf The Logic of Provability\rm , Cambridge University Press, Cambridge (1993). 
%DOI: \url{http://dx.doi.org/10.1017/CBO9780511625183}

\bibitem{bulseg} R.~A.~Bull and K.~Segerberg, {\em Basic modal logic}, pp. 1--88 in \cite{gab} and also pp. 1--82 in \cite{gab2}. 
%\url{http://dx.doi.org/10.1007/978-94-009-6259-0_1}

\bibitem{cz} A.~Chagrov. and M.~Zakharyaschev, \bf Modal Logic\rm , Clarendon Press, Oxford (1997).

\bibitem{fitting} M.~Fitting, \bf Proof Methods for Modal and Intuitionistic Logics\rm , D. Reidel, Dordrecht (1983).

\bibitem{gab} D.~Gabbay and F.~Guenthner (eds.), \bf Handbook of Philosophical Logic, vol. II: 
Extensions of Classical Logic\rm , D. Reidel, Dordrecht (1984).
%\url{http://dx.doi.org/10.1007/978-94-009-6259-0}

\bibitem{gab2} D.~Gabbay and F.~Guenthner (eds.), \bf Handbook of Philosophical Logic, vol. 3, 2nd ed.\rm , Springer, Dordrecht 2001. 

\bibitem{ham-lo} J.~D.~Hamkins and B.~L\"{o}we, {\em The modal logic of forcing}, \bf Transactions of the American Mathematical Society\rm ,
360 (2007), pp. 1793--1817.

\bibitem{hc1} G.~E.~Hughes and M.~J.~Cresswell, \bf An Introduction to Modal Logic\rm , Methuen, London (1968).

\bibitem{hc} G.~E.~Hughes and M.~J.~Cresswell, \bf A Companion to Modal Logic\rm , Methuen, London (1984).

\bibitem{humber} L.~Humberstone, \bf Philosophical Applications of Modal Logic\rm , College Publications, London (2015).

\bibitem{inoue1} T.~Inou\'{e}, {\em Partial interpretation of Le\'{s}niewski's epsilon in modal and intensional logics}, (abstract), 
\bf The Bulletin of Symbolic Logic\rm , 1 (1995), pp. 95--96. (I decided not to publish the full paper of this abstract, 
because \cite{blass} has been published and the essence of it is contained in \cite{inoue3}.) 
%DOI: 10.2307/420948 %\url{https://www.jstor.org/stable/420948}

\bibitem{inoue3} T.~Inou\'{e}, {\em Partial interpretations of Le\'{s}niewski's epsilon in von Wright-type deontic logics and provability logics}, 
\bf Bulletin of the Section of Logic\rm , 24 (1995), pp. 223--233.

\bibitem{inoue-Blass} T.~Inou\'{e}, \textit{On Blass translation for Le\'{s}niewski's propositional ontology and modal logics}, \bf Studia Logica\rm , (2021), https://doi.org/10.1007/s11225-021-09962-1. Also, arXiv:2006.15421v2 [math.LO], October 31, 2020\rm . The paper version will appear soon.
%DOI: \url{https://doi.org/10.1007/s11225-021-09962-16}

\bibitem{ishi} A.~Ishimoto, {\em A propositional fragment of Le\'{s}niewski's ontology}, \bf Studia Logica\rm , 36 (1977), pp. 285--299. 
%DOI: \url{http://dx.doi.org/10.1007/BF02120666}

\bibitem{iwa} B.~Iwanu\'{s}, {\em On Le\'{s}niewski's elementary ontology}, \bf Studia Logica\rm , 31 (1972), pp. 73--125.
%DOI: \url{http://dx.doi.org/10.1007/BF02120531}

\bibitem{koishi} M.~Kobayashi and A.~Ishimoto, {\em A propositional fragment of Le\'{s}niewski's ontology and its formulation by the tableau method}, 
\bf Studia Logica\rm , 41 (1982), pp. 181--195.
%DOI: \url{https://doi.org/10.1007/BF00370344}

\bibitem{ono1} H.~Ono, \bf Proof Theory and Algebra in Logic\rm, Springer, Singapore (2019).

\bibitem{pog} F.~Poggiolesi, \bf Gentzen Calculi for Modal Propositional Logic\rm, Springer, Dordrecht (2011).

\bibitem{savateev2021} Y.~Savateev and D.~Shamkanov,  {\em Non-well-founded proofs for the 
Grzegorczyk modal logic}, \bf The Review of Symbolic Logic\rm, 14 (2021), pp. 22--50.
%DOI: \url{https://doi.org/10.1017/S1755020319000510}

\bibitem{slu} J.~S\l upecki, {\em S. Le\'{s}niewski's calculus of names}, 
\bf Studia Logica\rm , 3 (1955), pp. 7--71.
%DOI: \url{http://dx.doi.org/10.1007/BF02067245}.

\bibitem{takano0} M.~Takano, {\em A semantical investigation into Le\'{s}niewski's axiom of his ontology}, 
\bf Studia Logica\rm , 44 (1985), pp. 71--77.
%DOI: \url{https://doi.org/10.1007/BF00370810}

\bibitem{urbaniak-book} R.~Urbaniak, \bf Le\'{s}niewski's Systems of Logic and Foundations of Mathematics\rm , Springer, Cham (2014).
%\url{http://dx.doi.org/10.1007/978-3-319-00482-2}

\end{thebibliography}
\end{document}